\documentclass{amsart}
\usepackage{amsmath}
\usepackage{amsfonts}
\usepackage{amssymb}
\usepackage{amscd}
\usepackage{graphics}
\usepackage{graphicx}
\usepackage{xypic}
 \input{xy}
  \xyoption{all}
\usepackage[lite]{amsrefs}

\usepackage{algorithmic}
\usepackage{clrscode} 

\newtheorem{theorem}{Theorem}

\newtheorem{lemma}[theorem]{Lemma}
 
\newtheorem{conjecture}{Conjecture}
\newtheorem*{theorem*}{Theorem}
\newtheorem*{question}{Question~Q(\textit{G,A,$\boldsymbol \ell$})}
\newtheorem*{alg}{Algorithm}

\theoremstyle{definition}
\newtheorem{definition}[theorem]{Definition}

\newtheorem{remark}[theorem]{Remark}

\title{Unstable Analogues of the Lichtenbaum-Quillen Conjecture}
\author{Marian Anton}
\address[Marian Anton]{CCSU\\
1516 Stanley St\\
New Britain, CT 06050\\
and IMAR\\
P.O. Box 1-764\\
014700\\
Bucharest, Romania}
\email{marianfanton9@gmail.com}

\author{Joshua Roberts}
\address[Joshua Roberts]{Piedmont College\\
165 Central Ave\\
Demorest, GA  30535}
\email{jroberts@piedmont.edu}
\urladdr{http://www.piedmont.edu/math/jroberts}

\begin{document}

\maketitle

\section{Introduction}

This survey is mostly concerned with unstable analogues
of the Lichtenbaum-Quillen Conjecture. The Lichtenbaum-Quillen
conjecture (now implied by the Voevodsky-Rost 
Theorem \cites{voevodsky:11}) attempts to describe the 
algebraic $K$-theory of rings of integers in number 
fields in terms of much more accessible ``\'{e}tale
models''. Suitable versions of the conjecture 
predict the cohomology of infinite general linear
groups of rings of $S$-integers over suitable number
fields; our survey focuses on an unstable version
of this form of the conjecture. This survey is not 
written for experts in algebraic $K$-theory and should
be seen as a computational viewpoint into the 
mysteries of the cohomology of $S$-arithmetic groups. 

\section{The General Problem}

The context of our survey can be described as follows:
let $G$ be a group scheme, $A$ a ring and $\ell$ a 
prime number. In the 1980's  
Dwyer and Friedlander \cites{dwyer:85}
introduced ``\'{e}tale models'' $BG^{et,\ell}(A)$ for 
the classifying space $BG(A)$ which come with natural
maps 
$$
f = f_{G,A,\ell}:BG(A) \to BG^{et,\ell}(A).
$$

These models have turned out to be particularly 
powerful in the context of algebraic $K$-theory, i.e.,
in the case of $G=GL_\infty$, in which the maps
$f_{G,A,\ell}$ are supposed to induce an isomorphism
in mod-$\ell$ cohomology for suitable primes $\ell$ and
certain rings $A$ of interest in number theory. In 
particular, this holds if $A=\mathbb Z[1/2]$ and
$\ell=2$ as implied by \cites{rognes:00}. 
In other words, the following question has
an affirmative answer in the case that $G=GL_\infty, 
A=\mathbb Z[1/2]$ and $\ell=2$.

\begin{question}
Does $f_{G,A,\ell}$ induce an isomorphism in mod-$\ell$
cohomology?
\end{question}

The situation turns out to be less favorable in the 
case of $G=GL_n$ and $G=O_n$. While the answer to question
 \textbf{Q(\textit{G,A,$\boldsymbol \ell$})} is 
affirmative if $G=GL_n$, $A=\mathbb Z[1/2]$ and $\ell=2$ as long
as $n \le 3$ as proven in \cites{henn:99, mitchell:92}, 
it was shown by Dwyer \cites{dwyer:98} that the answer 
is negative as soon as $n \ge 32$. In an unpublished
work of Lannes and Henn dating from the late 1990's
they showed that for $G=O_n, A=\mathbb Z[1/2]$ and 
$\ell=2$ the answer is yes if and only if $n \le 14$,
and for $G=GL_n, A=\mathbb Z[1/2]$ and $\ell=2$
the answer is no as soon as $n \ge 14$.

In our publications, we study question
\textbf{Q(\textit{G,A,$\boldsymbol \ell$})} in the case of $G=GL_n$, $\ell$ an odd
regular prime in the sense of number theory
and $A=\mathbb Z[1/\ell, \zeta_\ell]$, where $\zeta_\ell$
denotes a primitive $\ell$-root of unity. It is worthwhile
to note that the Voevodsky-Rost Theorem implies that the
question in this case has a positive answer if $G=GL_\infty$.

In \cites{anton:99} we study the case of the prime $3$. We show
by explicit calculation of 
$H^*(GL_2(\mathbb Z[1/3,\zeta_3]);\mathbb F_3)$ that 
the answer is 
positive if $\ell = 3$ and $n=2$. 

\begin{theorem}
$H^*(GL_2(\mathbb Z[1/3, \zeta_3]);\mathbb F_3) \approx \mathbb F_3[c_1, c_2] \otimes \Lambda_{\mathbb F_3}(e_1,e_2)\otimes \Lambda_{\mathbb F_3}(e_1', e_2')$.
\end{theorem}

Here $\Lambda_{\mathbb F_3}(e_1, e_2)$ denotes the exterior
algebra generated over $\mathbb F_3$ by the elements $e_1$ and 
$e_2$ where the degree of $e_j$ is $2j-1$ and similarly for 
$e_1'$ and $e_2'$. The polynomial algebra is generated by 
$c_1$ and $c_2$ where the degree of $c_j$ is $2j$.

This is the first 
non-trivial positive answer to the question in 
the case that $\ell$ is an odd prime. This calculation
requires considerable technical strength to combine
efficiently existing methods such as the spectral 
sequence associated to the action of $SL_2(\mathbb Z[\zeta_3])$ on a suitable $2-$dimensional contractible
complex as well as other spectral sequences 
converging to the high dimensional cohomology.
In the 
same paper we use the strategy used by Dwyer in \cites{dwyer:98}
 to show that the answer to 
\textbf{Q($\boldsymbol{GL_n,\mathbb Z[1/3,\zeta_3],3}$)} is negative as soon 
as $n \ge 27$. More precisely,

\begin{theorem}
The following homomorphism
$$
f^* \colon H^*(BGL_n^{et, \ell}(\mathbb Z[1/3, \zeta_3]);\mathbb F_3) \to H^*(BGL_n(\mathbb Z[1/3,\zeta_3]);\mathbb F_3)
$$
is injective for all $n \ge 0$ but not surjective for $n \ge 27$.
\end{theorem}

  As shown in \cites{henn:95} this disproves at 
the same time a conjecture of Quillen \cites{quillen:71} which 
predicted that the group cohomology $H^*(GL_n(\mathbb Z[1/3, \zeta_3]);\mathbb F_3)$ is a free module over 
$H^*(BGL_n(\mathbb C); \mathbb F_3)$ if the module
structure is induced from an embedding of the 
ring $\mathbb Z[1/3, \zeta_3]$ into the complex numbers
$\mathbb C$. Here $BGL_n(\mathbb C)$ is the classifying 
space of the Lie group $GL_n(\mathbb C)$ and not of the 
discrete group.

In \cites{anton:01} we show that Dwyer's strategy of \cites{dwyer:98}
can be applied to show that the answer to 
\textbf{Q(\textit{G,A,$\boldsymbol \ell$})} is negative
if $G=GL_n$ for $n$ sufficiently large if $\ell$ is 
an odd regular prime and if $A=\mathbb Z[1/\ell, \zeta_{\ell}]$. 

\begin{theorem}\label{reg-prime}
If $\ell$ is an odd regular prime, then the map
$$
f = f_n \colon BGL_n(\mathbb Z[1/\ell, \zeta_\ell]) \to BGL_n^{et, \ell}(\mathbb Z[1/\ell, \zeta_\ell])
$$
does not induce an isomorphism on mod-$\ell$ cohomology 
for $n$ sufficiently large.
\end{theorem}

The proof requires a careful analysis
of the \'{e}tale model and uses (as in \cites{dwyer:98}) deep
results from the homotopy theory of classifying spaces,
which ultimately depend on the solution of the 
Sullivan conjecture.

\begin{remark}
The map $f_n$ in Theorem \ref{reg-prime} induces an injection
on mod $\ell$ cohomology for all $n \ge 0$ by \cites{dwyer:94}
but not a surjection for $n \ge N_\ell$ where $N_\ell$ depends
on $\ell$.
\end{remark}

In \cites{anton:03} the prime $\ell$ and the ring $A$ are as 
in \cites{anton:01}. The focus is on analyzing the image $I_n$
of $H^*(BGL_n(A); \mathbb F_\ell)$ in $H^*(BD_n(A);\mathbb F_\ell)$ with respect to 
the restriction homomorphism $res_n$ associated 
to the inclusion of the subgroup $D_n(A)$ of diagonal
matrices into $GL_n(A)$. This is compared with the 
image $M_n$ of the composition of $res_n \circ f_n^*$.
It is clear that $M_n \subset I_n$. Furthermore, 
the \'{e}tale model can be explicitly described
and $M_n$ can be explicitly calculated as in \cites{dwyer:94}. 

\begin{theorem}
If $\ell$ is an odd regular prime, then $res_n \circ f_n^*$ induces an isomorphism
$$
H^*(BGL^{et, \ell}_n (\mathbb Z[1/\ell, \zeta_l]); \mathbb F_\ell) \approx M_n = \mathbb F_\ell[c_1, c_2, \dots,
c_r] \otimes \Lambda_{\mathbb F_\ell} (e_{i,1}, e_{i,2}, \dots, e_{i,r})^{\otimes^r_{i=1}}
$$
where $r=\frac{1}{2}(\ell + 1)$, the generator $c_j$ has degree $2j$ and $e_{i,j}$ has degree $2j-1$.
\end{theorem}

The main 
result of the paper says that $M_2=I_2$ implies
$M_n=I_n$ for all $n \ge 2$; in other words, if one 
can verify for $n=2$ that the image of the restriction
is not larger than the image of the \'{e}tale model,
then this holds true for every $n \ge 2$. 

\begin{theorem}\label{invariant}
If $\ell$ is an odd regular prime, then the images of the maps
$res_n \circ f_n^*$ and $res_n$ agree i.e., $M_n = I_n$ for
all $n \ge 0$, if and only if $M_2 = I_2$.
\end{theorem}

The proof
of this statement relies on the solution to a 
problem in invariant theory involving the canonical
(signed) permutation action of the symmetric group
$\Sigma_n$ on graded rings of the form 
$\left(\mathbb F_\ell[x] \otimes \Lambda_{\mathbb F_\ell}(y_1, \dots, y_r)\right)^{\otimes n}$.

\section{Homological Symbols}

Unfortunately the verification of $M_2=I_2$ is 
extremely hard to verify as soon as the prime $\ell$ is 
bigger than 3. In order to pin down where the difficulties lie, in 
\cites{anton:08} we concentrate our attention on the image 
of the canonical homomorphism 
$$
H_i(BD_j(\mathbb Z[1/\ell, \zeta_\ell]); \mathbb F_\ell) \to H_i(BGL_j(\mathbb Z[1/\ell, \zeta_\ell]); \mathbb F_\ell)
$$
induced by the inclusion of the subgroup of diagonal
matrices $D_j \subset GL_j$. As $i$ and $j$ vary, one gets a 
homomorphism of bi-graded algebras.

\begin{definition}We call the image of this homomorphism the \textit{algebra of 
homological symbols} and denote it by $KH_{**}(\mathbb Z[\zeta_\ell, 1/\ell])$. 
\end{definition}

We proposed to compare this 
to an \textit{algebra of \'{e}tale homological 
symbols} similarly defined in terms of the \'{e}tale
unstable models for $BGL_j(\mathbb Z[\zeta_\ell, 1/\ell])$ and we denote that algebra by 
$KH^{et}_{**}(\mathbb Z[\zeta_\ell, 1/\ell])$. We had 
explicitly determined the algebra $KH^{et}_{**}$ in 
terms of generators and relations as follows.

The mod $\ell$ group homology of $GL_1$ is the homology of the bar resolution $(B_*,\partial)$
where $B_i$ is the vector space over $\mathbb F_\ell$ generated by the symbols $[x_1|\dots|x_i]$ with
$x_i \in GL_1\setminus \{1\}$ and $\partial$ is the boundary homomorphism
$$
\partial [x_1|\dots|x_i]=[x_2|\dots|x_i] + \displaystyle \sum_{j=1}^{i-1}(-1)^j [x_1|\dots|x_j x_{j+1}| \dots 
|x_i] + (-1)^i[x_1|\dots|x_{i-1}]
$$
where we delete all the terms with $x_jx_{j+1}=1$. We define the shuffle product of two symbols by the 
formula
$$
[x_1|\dots|x_i] \wedge [x_{i+1}|\dots|x_{i+s}] = \displaystyle \sum_\sigma 
(-1)^\sigma[x_{\sigma(1)}|\dots|x_{\sigma(i+s)}]
$$
where $\sigma$ runs over all the permutations of $\{1, \dots, i+s\}$ that shuffle $\{1, \dots, i\}$
with $\{i+1, \dots, i+s\}$ and $(-1)^\sigma$ denotes the signature of $\sigma$.

By \cites{washington:87}, $GL_1$ is the Abelian group generated by the set of cyclotomic units
$$
S=\{-\zeta_\ell, 1-\zeta_\ell, 1-\zeta_\ell^2, \dots, 1-\zeta_\ell^r\}
$$
subject to one relation $(-\zeta_\ell)^{2\ell}=1$. Here we note that $\ell=2r+1$. By a K\"{u}nneth isomorphism, the bi-graded algebra 
$H_i(BD_j(\mathbb Z[\zeta_\ell, 1/\ell]); \mathbb F_\ell)$ is the bi-graded tensor algebra generated by the 
cycles
\begin{equation}\label{eq}
[\zeta_\ell]^{(s)} \wedge \langle v_1, \dots, v_i \rangle := \displaystyle \sum^{\ell -1}_{i_1, \dots, i_s=1}
[\zeta_\ell^{i_1} | \zeta_\ell | \dots | \zeta_\ell^{i_s} | \zeta_\ell] \wedge [v_1] \wedge \dots \wedge [v_i]
\end{equation}
where $\{v_1, \dots, v_i\} \subset S$ is any subset of cyclotomic units, $s$ is any nonnegative integer, 
and the bi-degree of this cycle is $(i+2s, 1)$.

\begin{definition}
A cycle of the form (\ref{eq}) is an \textit{\'{e}tale obstruction cycle} if $i-s$ is $ > 0$ and even and it is 
an \textit{odd cycle} if $i+s$ is odd.
\end{definition}

\begin{theorem}\label{algebras}
If $\ell$ is an odd regular prime, then $KH^{et}_{**}$ is the bi-graded algebra generated by cycles of the 
form (\ref{eq}) and subject to the relations
$$
\rho_*(t_*(z) \otimes z')=0, \; t(u) = u^{-1} \times u, \; \rho( u \times v \times w)=uw \times vw
$$
where $z$ runs over the \'{e}tale obstruction cycles and $z'$ over the odd cycles.
\end{theorem}
Here $t_*$ and $\rho_*$ are induced at the level of cycles by the group homomorphisms 
$t \colon GL_1 \to GL_1^{\times 2}$ and $\rho \colon GL_1^{\times 3} \to GL_1^{\times 2}$ defined in the theorem. Moreover there is a canonical surjection of bi-graded algebras
$$
KH_{**}(\mathbb Z[\zeta_\ell, 1/\ell]) \to KH^{et}_{**}(\mathbb Z[\zeta_\ell, 1/\ell])$$
and we conjecture 

\begin{conjecture}\label{conj}
This homomorphism is an isomorphism for $\ell$ a regular odd prime.
\end{conjecture}

We verify this conjecture in the case of the prime
$\ell =3$ and indicate a program how one might 
attack the case of larger primes. The case of 
each prime would have to be dealt with separately
and the computational complexity would grow
rapidly with the prime. As a consequence of Theorem \ref{algebras} we have the following result \cites{anton:03}:

\begin{theorem}\label{etale} The Conjecture \ref{conj} is true if and only if $t_*(z)$ is null-homologous in 
$H_*(SL_2(\mathbb Z[1/\ell, \zeta_\ell]); \mathbb F_\ell)$ for all \'{e}tale obstruction cycles $z$.
\end{theorem}

\begin{remark}
 The Conjecture \ref{conj} is equivalent to $M_2=I_2$ as in Theorem \ref{invariant}. Observe that the image of the homomorphism $t$ is contained in the diagonal subgroup of $SL_2$ and there are finitely many {\'e}tale obstruction cycles to check. 
\end{remark}

\section{Some algorithms}

These verifications are usually not accessible 
without the algorithms developed by the 
second author in \cites{roberts:10}. Namely, given $G$ a finitely-presented group and $k$ a
finite field, the paper exploits a formula due to Hopf to algorithmically
give upper bounds on the dimension $d$ of the second homology group of $G$ as a vector space over $k$.

In particular, if $G = F/R$ is a finite presentation, we have the short exact sequence
$$
1 \to [F,R] \to R \cap [F,F] \to H_2(G) \to 1.
$$
Through a series of short exact sequences, we show that the dimension of $H_2(G,k)$ is less than 
or equal to $a+b-c+e$, where $a$ is the dimension of 
$\mathrm{Tor}(H_1(G),k)$, $\ell^b$ and $\ell^c$ are the 
orders of the $\ell$-primary subgroups of $F/R[F,F]$ and 
$F/R^\ell[F,F]$ respectively, $e$ is the vector space 
dimension of $k \otimes R/[F,R]$, and $\ell$ is the
characteristic of $k$. All of these numbers except for $e$ can be calculated 
by algorithms given in \cites{roberts:10}. In general, $e$ can only be estimated 
from above using the following series of algorithms given in pseudocode.

\begin{alg}
$\proc{ReduceWord}(F, R, Z, R', \ell)$
\begin{algorithmic}[1]
\REQUIRE Free Group $F$, Relators $R$, Test Word $z$, Sublist $R'$ of $R$, Prime $\ell$
\ENSURE Reduced word of $z$ in $F/[F,R]R^\ell R'$
\STATE $G:=F/[F,R]R^\ell R'$
\STATE $RG:=$Rewriting system for $G$
\STATE $x:=$Reduced word of $(z)$ in the rewriting system $RG$
\RETURN $x$
\end{algorithmic}
\end{alg}

This algorithm attempts to reduce a given test word in a finitely-presented group using a rewriting system.
We utilized the KBMAG package implemented on GAP for this purpose.

\begin{alg}
$\proc{FindBasis}(F,R,\ell,R')$
\begin{algorithmic}[1]
\REQUIRE Free Group $F$, Relators $R$, Prime $\ell$, Sublist $R'$ of $R$
\ENSURE Size of a generating set for $[F,R]R^\ell R'/[F,R]R^\ell$
\STATE $X:=R'$
\FOR{$x \in X$}
\STATE $x':=$\textsc{ReduceWord}$(F,R,x,\text{Difference}(X,[x]),\ell)$ \COMMENT{Difference$(A,B)$ is the complement of $B$ in $A$}
\IF{$x'=$ identity}
\STATE $X:=$Difference$(X,[x])$
\ENDIF
\ENDFOR
\RETURN Size$(X)$
\end{algorithmic}
\end{alg}

The algorithm tries to determine the linear independence of elements $x$ in a generating set $R'$ with respect to $R'-\{x\}$ in $[F,R]R^\ell R'/[F,R]R^\ell$. If $x$ is determined to be 
dependent on $R'-\{x\}$, it is removed from $R'$. The end result will be a list of generators which 
are, at least potentially, linearly independent.

\begin{alg}
$\proc{SecondHomologyCoefficients}(F, R, \ell, R')$ \label{masteralg}
\begin{algorithmic}[1]
\REQUIRE Free Group $F$, Relators $R$, Prime $\ell$, Sublist $R'$ of $R$ generating $R/[F,R]R^\ell$
\ENSURE An integer $d$ such that dim $\left(H_2(G;k)\right) \le d$
\STATE $a:=$ \textsc{Tor}$(F,R,\ell)$
\STATE $b:=$ \textsc{PrimePrimaryRank}$(F,R[F,F],\ell)$
\STATE $c:=$ \textsc{PrimePrimaryRank}$(F,R^p[F,F],\ell)$
\STATE $e:=$ \textsc{FindBasis}$(F,R,\ell,R')$
\STATE $d:=a+b-c+e$
\RETURN $d$
\end{algorithmic}
\end{alg}

Based on finite presentations for groups $G = SL_2(\mathbb Z[1/\ell, \zeta_\ell])$ in 
\cites{anton:08} and our algorithms in \cites{roberts:10} we prove that $d = 0$ for $\ell = 3, 5$ and $d \le 6$ for $\ell = 7$. 

\begin{theorem}\label{sl2}
The homology group $H_2(SL_2(\mathbb Z[1/\ell,\zeta_\ell]);\mathbb F_\ell)=0$ if $\ell=3,5$ and it is at most $6$-dimensional as a vector space over $\mathbb F_\ell$ if $\ell=7$.
\end{theorem} 

To prove this theorem we define a group $SE_2$ for $\ell=2r+1$ by a finite presentation
$$1\to R\to F\to SE_2\to 1$$
where $F$ is the free group given by the generators $z,u_1,...,u_r,a,b$ and $R$ is the normal subgroup given by the following relations
$$z^\ell=[z,u_i]=[u_i,u_j]=a^4=[a^2,z]=[a^2,u_i]=[b_s,b_t]=c(I)^3=1$$ 
$$a=zaz=u_i au_i, b^3=a^2=b_0b_1...b_{2r}, b_t ^\ell=w^{-1}b_t^{(-1)^r}w,$$
$$ba^2=u_ibz^{-ri}b^{-1}b_0^{-1}z^{ri}bz^{-i}u_i$$
where $i,j\in\{1,2,...,r\}$, $s,t\in\{0,1,2,...,2r\}$, and $I\subset\{1,2,...,r\}$ nonempty. The notations are as follows
$$b_t:=z^{rt}bz^{rt}a, w:=z^cu_1u_2...u_r, c(I):=(\prod_{t=0}^{2r}b_t^{c_t(I)})a^{-1}\prod_{i\in I}u_i$$ 
where $c,c_t(I)\in\mathbb Z$ with $c\ge0$ minimal such that
$$
2c\equiv r^2+\frac{r(r+1)}{2}\mod\ell, \prod_{i\in I}(1-\zeta_\ell^i)=\sum_{t=0}^{2r}c_t(I)\zeta_\ell^t.
$$
The conventions are $b_t=b_s$ if $s\equiv t\mod\ell$ and $[x,y]=xyx^{-1}y^{-1}$. Then there is a group homomorphism $SE_2\to SL_2$ such that Conjecture \ref{conj} is true if the cycles
$$
[z]^{(s)}\wedge[e_1]\wedge...\wedge[e_i]
$$
are null-homologous in $H_{3s+2j}(SE_2;\mathbb F_\ell)$ for each pair of nonnegative integers $(s,j)$ and subset $\{e_1,...,e_i\}\subset\{z,u_1,...,u_r\}$ with $i=s+2j$ and $j>0$. Our algorithms apply in the case $(s,j)=(0,1)$ as follows

\begin{lemma} The cycle $[e_1]\wedge[e_2]$ is null-homologous in $H_2 (SE_2;\mathbb F_\ell)$ if and only if $[e_1,e_2]\in [F,R]R^\ell$.\end{lemma}

Combining  Theorems \ref{etale} and \ref{sl2} we conclude that Conjecture \ref{conj} is true if $\ell=3$. For $\ell\ge 5$ our algorithms need to be improved to handle higher computational complexity for our conjecture. Nevertheless, there are a couple of new applications we consider in our work in progress \cites{antonroberts:12}. 

We are constructing, as a by product of the calculations relevant for this paper, a database for low dimensional
group homology of linear groups. This will be extended to other finitely-presented groups of interest
in number theory and computational topology. An initial set of these calculations is given in Table 1.
 While the results in the table are not new, they were previously found by a variety of methods, many
of which are not computational. A $\le$ indicates that only 
an upper bound was found.

\begin{table}[h]
\centering
$$
\begin{array}{|l|c|c|c|c|}
\hline
                   & H_2(-;\mathbb F_2) & H_2(-;\mathbb F_3) & H_2(-;\mathbb F_5) & H_2(-;\mathbb F_7) \\\hline
GL_2(\mathbb Z)           &    \le 4        & \le 2            &  \le 2         &   \le 2           \\\hline
SL_2(\mathbb Z)           &     \le 2       &   \le 2          &  \le 1         &   \le 1           \\\hline
SL_2(\mathbb Z_2)         &    1       &     0        &   0         &    0       \\\hline
SL_2(\mathbb Z_3)         &    0        &    1        &   0         &    0       \\\hline
SL_2(\mathbb Z_5)         &    0        &    0         &  1         &    0       \\\hline
SL_2(\mathbb Z[i])        &    1        &   0          &   0        &    0          \\\hline
SL_2(\mathbb Z[\omega]), \omega^3=-1   &     \le 1       &  \le 2           &    \le 1       &   \le 1          \\\hline
SL_2(\mathbb Z[\sqrt{-5}])&    \le 3        &    \le 3         &    0       &      0        \\\hline
PSL_2(\mathbb Z)          &    \le 1        &    \le 1         &    0       &      0        \\\hline
\end{array}
$$

\caption{Dimensions of Second Homology Groups}
\end{table}

Also in \cites{antonroberts:12} we explain how the 
algorithms in \cites{roberts:10} can be adapted 
to find generators of $H_2$. For example, results 
in \cites{anton:08} show that $SL_2(\mathbb Z[1/7, \zeta_7])$ 
has presentation with generators  $\{z,u_1,u_2,u_3,a,b,b_0,b_1,b_2,b_3,b_4,b_5,b_6,w\}$
and a set of 64 relators.

\begin{theorem}
As a vector space over $\mathbb F_7$, $H_2(SL_2(\mathbb Z[1/7, \zeta_7]); \mathbb F_7)$ is generated by
the following set of relators
$$
\begin{array}{l}
\{ z a z a^{-1}, \\
u_1 u_2 u_1^{-1} u_2^{-1},\\
u_1 a u_1 a^{-1},\\
 b_1 a^{-1} b_1 a^{-1} b_1 a,\\
  u_2 b_1 z^{-2} b_1 z b_1^{-2} a^{-1} u_2 z^{-1} a^{-1} b_1^{-1}, \\
  u_3 b_1 z^2 b_1^{-2} z^2 b_1 a^{-1} z^{-2} u_3 z^{-1} a^{-1} b_1^{-1} \}.
\end{array}
$$
\end{theorem}

To summarize, significant contributions were made,
which help clarify the limitation of \'{e}tale
models for finite general linear groups over rings 
of $S$-integers of number fields, but much remains
to be done in low dimensional group homology.

\end{document}